\documentclass[11pt, reqno, a4paper]{amsart}

\usepackage{amssymb, amsmath, amscd}
\usepackage{amsthm}

\usepackage{verbatim}

\usepackage{color}

\usepackage{etex}
\usepackage{dsfont}
\usepackage[matrix, arrow, curve, frame, arc]{xy}
\usepackage{pgf}
\usepackage{tikz-cd}
\usetikzlibrary{arrows,shapes,automata,backgrounds,petri,calc}
\usepackage[square,sort,comma,numbers]{natbib}
 \usepackage{url}
 \usepackage{hyperref}
 \usepackage[shortlabels]{enumitem} 

\usepackage{adjustbox}

\usepackage{algorithm,algpseudocode}

\newcommand{\tikzAngleOfLine}{\tikz@AngleOfLine}
\def\tikz@AngleOfLine(#1)(#2)#3{%
\pgfmathanglebetweenpoints{%
\pgfpointanchor{#1}{center}}{%
\pgfpointanchor{#2}{center}}
\pgfmathsetmacro{#3}{\pgfmathresult}%
}

\newcommand{\bS}{\mathbb{S}}

\newcommand{\bN}{\mathbb{N}}
\newcommand{\bM}{\mathbb{M}} 
\newcommand{\bfQ}{{\bf Q}}
\newcommand{\bbQ}{\mathbb{Q}}
\newcommand{\bZ}{\mathbb{Z}}

\newcommand{\ba}{\bar{\alpha}}
\newcommand{\La}{\Lambda}

\newcommand{\cA}{\mathcal{A}}
\newcommand{\cS}{\mathcal{S}}

\newcommand{\Adj}{\operatorname{Ad}}
\newcommand{\Arr}{\operatorname{Ar}} 
\newcommand{\op}{\operatorname{op}}

\begin{document}

\newtheorem{defi}{Definition}[section] 
\newtheorem*{deff}{Definition}
\newtheorem{rem}[defi]{Remark}
\newtheorem{prop}[defi]{Proposition}
\newtheorem{ques}[defi]{Question}
\newtheorem{lemma}[defi]{Lemma}
\newtheorem{cor}[defi]{Corollary}
\newtheorem{thm}[defi]{Theorem}
\newtheorem{expl}[defi]{Example} 

\newcounter{mac}

\parindent0pt

\title[]{Periodicity shadows II. Computational aspects$^*$ \footnote{\tiny$^*$This research was financed from the grant no. 
2023/51/D/ST1/01214 of the Polish National Science Center}}

\author[J. Bia{\l}kowski]{Jerzy Bia{\l}kowski}
\address[Jerzy Bia{\l}kowski]{Faculty of Mathematics and Computer Science, 
Nicolaus Copernicus University, Chopina 12/18, 87-100 Torun, Poland}
\email{jb@mat.umk.pl} 

\author[A. Skowyrski]{Adam Skowyrski}
\address[Adam Skowyrski]{Faculty of Mathematics and Computer Science, 
Nicolaus Copernicus University, Chopina 12/18, 87-100 Torun, Poland}
\email{skowyr@mat.umk.pl} 


\begin{abstract} This article provides the second part of the research initiated in \cite{Sko}, where we introduced 
and investigated so called periodicity shadows, which are special skew-symmetric matrices related to symmetric algebras 
with periodic simple modules. In \cite{Sko} we focused on theoretical aspects, whereas here we present complementary 
cosiderations concerning computational issues. Namely, we discuss an algorithm, which computes all tame periodicity 
shadows of given size (see Section \ref{sec:2}), and then present the lists of tame periodicity shadows of small sizes 
$n\leqslant 6$. \end{abstract}

\maketitle

\section{Introduction}\label{sec:1} 

By an algebra, we mean an indecomposable basic finite-dimensional $K$-algebra over an algebraically closed field 
$K$. It is known that $\Lambda$ has presentation of the form $\Lambda \cong KQ/I$, where $KQ$ is the path algebra 
of a quiver $Q$, and $I$ is so called admissible ideal. The quiver $Q$ is uniquely determined and it is called 
the {\it Gabriel quiver} of $\Lambda$, and denoted by $Q_\La$. All modules over an algebra $\La$ are assumed 
finitely generated, or equivalently, finite-dimensional. \medskip 

Tools invented in \cite{Sko} are designed to study special classes of algebras, i.e. symmetric periodic algebras 
of period four, or more generally, so called algebras of generalized quaternion type \cite{AGQT}, which constitute 
important subclasses of the class of selfinjective algebras. Except one place, we will not deal with algebras 
directly, focusing mainly on the combinatorial life of associated matrices. For precise definitions we refer to 
\cite[see Sections 1 and 2]{Sko} or \cite{Sk-Y}. \smallskip 

We recall that for any algebra $\La=KQ/I$ we have a complete set of pairwise orthogonal idempotents $e_i$, where 
$i$ runs through the vertex set $\{1,\dots,n\}$ of $Q$ and $e_i$ is the coset of the trivial path (of length $0$) 
at vertex $i$. It induces a decomposition of the identity $1_\Lambda=\sum_{i=1}^n e_i$ and there is an associated 
decomposition $\La=e_1\La\oplus\dots\oplus e_n\La$ of $\La$ into a direct sum of indecomposable projective 
$\La$-modules, where the modules $e_i\La$, for $i\in\{1,\dots,n\}$, form a complete set of pairwise non-isomorphic 
indecomposable projective $\La$-modules (see \cite{Sk-Y}). The Cartan matrix $C_\La$ of an algebra $\La$ is a $n\times n$ 
matrix with natural coefficients, whose columns are dimension vectors of projective modules $e_i\La$. It means 
that $C_\Lambda=[c_{ij}]$ with $c_{ij}=\dim_K e_j\La e_i$. Note that $C_\La$ is a symmetric matrix if the algebra 
$\La$ is symmetric. \medskip  

By a {\it periodic} module (of period $d$), we mean a (right) $\Lambda$-module $X$ which is periodic with 
respect to the syzygy operator $\Omega_\La$, that is, we have $\Omega_\La^d(X)\simeq X$ and $d$ is minimal number 
with this property. \medskip 

Let us now briefly introduce the notion of a periodicity shadow, for short, just shadow; more details in \cite{Sko}.  
Origins of this notion come from studying the properties of the adjacency matrices of Gabriel quivers 
defining symmetric algebras with all simple modules of period four. For a quiver $Q$, we denote by 
$\Arr_Q=[\alpha_{ij}]\in\bM_n(\bN)$ the {\it arrow matrix} of $Q$, where $\alpha_{ij}$ is the number of arrows 
$i\to j$ in $Q$. Then the {\it (signed) adjacency matrix} $\Adj_Q=[a_{ij}]$ of $Q$ is defined as 
$$\Adj_Q=\Arr_Q-\Arr_Q^T,$$ 
that is, $a_{ij}$ counts the difference between the number of arrows $i\to j$ and the number of arrows $j\to i$ 
in $Q$. The study of the (signed) adjacency matrices in the context of algebras was motivated by the property 
observed in \cite[Theorem 2.1]{Sko}, where we proved that for any symmetric algebra $\La=KQ/I$, $Q=Q_\La$, 
with all simple $\La$-modules periodic of period $4$, its Cartan matrix $C=C_\Lambda$ satisfies the following 
identity: 
$$\Adj_Q\cdot C=0.$$ \medskip 

Roughly speaking, the above property can be used as a tool to describe Gabriel quivers of tame symmetric algebras 
with simples periodic of period $4$. Namely, we may view it as a equation, where $\Adj_Q$ is a fixed skew-symmetric 
matrix (representing unknown Gabriel quiver), and $C$ is a matrix variable, and then look for all matrices $A$ such 
that $AC=0$ has solutions which are symmetric matrices with natural coefficients (and non-zero columns). We only 
mention that $A$ determines $Q$ up to loops and $2$-cycles, but there are precise rules which allow to reconstruct $Q$ 
from $A$ (see the Main Theorem of \cite{Sko}). In this paper, we only focus on the algorithm computing the shadows 
and we present complete lists for $n\leqslant 6$. \medskip 

Following \cite{Sko}, we recall the definition of a periodicity shadow. 

\begin{deff} \normalfont A matrix $A\in\bM_n(\bZ)$ is called a {\it periodicity shadow}, if the following holds: 
\begin{enumerate}
\item[PS1)] $A$ is a singular skew-symmetric matrix.  
\item[PS2)] $A$ does not admit a nonzero row containing integers of the same sign. 
\item[PS3)] There exists a symmetric matrix $C\in\bM_n(\bN)$ with non-zero columns, such that $AC=0$. 
\end{enumerate}\end{deff} 

A periodicity shadow $A$ is called {\it tame} if the following conditions hold: 

\begin{enumerate}
\item[T1)] all entries $a_{ij}$ of $A$ are in $\{-2,-1,0,1,2\}$. 
\item[T2)] each row of $A$ does not contain simultaneously: both $2$ and entry $\geqslant 1$ or both $-2$ and entry 
$\leqslant -1$. 
\item[T3)] each row of $A$ cannot have more than four $1$'s or more that four $-1$'s. 
\end{enumerate} 

We only mention that any tame symmetric algebra $\Lambda=KQ/I$ with all simple modules periodic of period $4$ 
gives rise to a tame periodicity shadow $\bS_\Lambda=\Adj_{Q_\La}$. \medskip 

To compute all tame periodicity shadows of a given size, we will actually compute a bit more, namely all matrices 
that are called shades. By a {\it shade}, we mean any integer matrix satisfying conditions PS1)-PS2) and T1)-T3). 
So a shade is a (tame) periodicity shadow if and only if it satisfies PS3). \medskip 

The algorithm described in the next section computes all basic shades of a given size $n$, where by basic we mean that 
we compute some subset of shades, from which every other shade can be obtained by permutation of rows and columns. 
We only mention that even for small number $n$ of vertices, say $n\leqslant 6$, it is a challenge, whereas for bigger 
$n$, it quickly becomes intractable in practice. All we know for $n\geqslant 8$ is that there exists a set of tame 
periodicity shadows, and there is an algorithm (based on recursive generating) which can compute it, at least in theory. 
The idea is to first generate the set $\cS(n)$ of all (basic) shades, and then compute all solutions of $AC=0$ with $C=C^T$, 
for any $A\in\cS(n)$, from which one can decide whether a given shade is a (tame) periodicity shadow. We will denote by 
$\bS(n)$ the list of all basic tame periodicity shadows $\bS(n)\subseteq\cS(n)$, obtained from $\cS(n)$ by removing 
matrices not satisfying PS3). \medskip 

For small $n$, we have $\bS(3)=\cS(3)$ and $\bS(4)=\cS(4)$, and these lists are presented in Section \ref{sec:3}. 
For $n\geqslant 5$ there are shades which are not tame periodicity shadows. Because of the length of the lists 
$\cS(5)$ and $\cS(6)$, we will restrict ourselves only to so called {\it essential shadows}, since any non-essential 
one is not a shadow of a tame algebra (see Section \ref{sec:2}). Full lists of essential shadows of sizes $n=5,6$ are 
presented in Sections \ref{sec:4}-\ref{sec:5}. For the reader interested in the complete lists $\cS(5)$ and 
$\cS(6)$ containing non-essential shadows and shades which are not shadows, we refer to the Appendix. \medskip 

Though the algorithm is useless in practice for $n\geqslant 9$, this tool was invented especially for small cases 
$n\leqslant 6$ to understand the structure of Gabriel quivers of small periodic algebras. This could be also an 
important conceptual tool to study general structure of such quivers, as one can see in \cite{Sko}. On the other hand, 
it is conjectured that for $n\geqslant 7$ the picture becomes more stable, and known examples of families of 
exceptional algebras appear only for $n=6$; see \cite{HTA,HSA,TT}. Last but not least, studying shades may be also 
important. However, shades which are not shadows (or essential shadows) do not come from an algebra, but understanding 
the difference at the combinatorial level could give new useful insights. \bigskip 

\section{The algorithm}\label{sec:2} 

In this section we will discuss technical issues around computation of tame periodicity shadows and present 
the algorithm (together with the proof of its correctness). \bigskip 

First naive method would be to search through the set $\bM^*_n$ of all skew-symmetric $n\times n$ matrices $A$ 
with $a_{ij}\in\{-2,-1,0,1,2\}$ (any such a matrix has $a_{ii}=0$), and check whether $A$ is a tame periodicity 
shadow, which requires verifying conditions T2)-T3) and PS1)-PS3). This can be done by hand, for $n=3$ and with 
a little bit more effort for $n=4$, but for higher $n$ it is unfeasible even with a help of a computer. \medskip 

The idea of the algorithm is to give a recursive procedure generating all the tame periodicity shadows (up to permutation 
of rows and columns and taking $-A=A^T$), instead of performing a naive search. Additionally, the output list is 
sorted according to some natural order. This gives an effective algorithm generating all tame periodicity shadows for 
$n\leqslant 6$. For $n=7,8$ computations are possible, but take a lot of time, whereas for $n\geqslant 9$, 
given the computational capabilities of today's computers (mainly due to the ``exploding'' numbers of both cases 
and results) it seems to be hopeless, or at least impractical. However, the tool was invented to handle especially 
small cases, and it is not our aim to compute all the shadows for bigger $n$. \smallskip 

Let $S_n\subset\bM^*_n$ be the subset of all tame periodicity shadows. We note that $S_n$ is closed under 
the action of permutations, by which we mean that given $A\in S_n$, we have $A'=P^TAP\in S_n$, for any 
permutation matrix $P$. We denote by $G$ the group of permutation matrices, and for any matrix $A$, we 
write $[A]$ for the $G$-orbit $\{P^TAP; \ P\in G\}$ of $A$ under this action. In other words, $S_n$ is a disjoint 
union of $G$-orbits, which are equivalence classes of congruence relation $\sim$, where $A \sim B$ if and 
only if $B=P^TAP$, for some $P\in G$. \smallskip 

We will frequently identify a tame periodicity shadow $A$ with the associated quiver $\bfQ_A$, which is the 
smallest quiver having $A$ as its signed adjacency matrix \cite[see Section 2]{Sko}. \smallskip 

Moreover, there is a natural isomorphism between the spaces of solutions of $AC=0$ and $BC=0$, if $A\sim B$ 
(preserving symmetricity and natural coefficients). In the language of associated quivers, it means that $\bfQ_A$ 
and $\bfQ_{B}$ are the same up to permutation of vertices encoded in $P$. Our aim is to find `basic' shapes, 
so we want to identify matrices in $S_n$, which belong to the same orbit of this action. Namely, we are interested in 
finding one particular set of representatives $\mathcal{A}=\{\mathbb{A}_1,\dots,\mathbb{A}_r\} \subset S_n$ such 
that $S_n$ is a disjoint union of orbits 
$$S_n=\bigcup_{i=1}^r[\mathbb{A}_i].$$ 
Given one $\cA$, any other subset with this property is determined up to permutation, that is, it consists of $r$ 
matrices of the form $P_i^T\mathbb{A}_iP_i$, for some permutation matrices $P_1,\dots,P_r$. \medskip 

Further, for any skew-symmetric integer matrix $A$, we have $A^T=-A$, so the nullspaces are equal: 
$\mathcal{N}(A^T)=\mathcal{N}(A)$. Hence, the matrices $A$ and $A^T$ induce the same set of solutions 
of matrix equation in PS3). Clearly, all the other conditions PS1)-PS2) and T1)-T3) are also preserved 
in $A^T=-A$, hence $A^T\in S_n$, if $A\in S_n$. Note also that taking $A^T$ is compatible with the 
action of permutations, i.e. $A^T \sim B^T$ for $A\sim B$, or equivalently, $[-A]=[-B]$, for $[A]=[B]$. 
By definition, at the level of quivers $\bfQ_{A^T}$ is the opposite quiver 
$$\bfQ_{A^T}=\bfQ_{-A}=\bfQ_A^{\op}$$ 
of $\bfQ_A$ (obtained from $\bfQ_A$ by reversing the orientation of all arrows). \medskip 

Now, using both of the above actions we can define what are the {\it basic} (tame) periodicity shadows. 
We start with the list $\bS=(\mathbb{A}_1,\dots,\mathbb{A}_r)$ of all representatives of $G$-orbits, 
and modify it in the following way: 

\begin{algorithmic}[H]
	\For{$j \mbox{ \bfseries from } 2 \mbox{ \bfseries to } r$}
	\If{$\exists_{i<j} [-\mathbb{A}_i]=[\mathbb{A}_j]$} 
	\State
	$S \gets S\setminus A_j$
	\EndIf   	  
	\EndFor
\end{algorithmic} \medskip

In this way, we remove from the representatives of $G$-orbits all the matrices $\mathbb{A}_j$ which are 
opposite to some earlier matrix $\mathbb{A}_i\in\mathcal{A}$, modulo permutation. The resulting list 
$\bS=(\bS_1,\dots,\bS_s)$ contains so called {\it basic} tame periodicity shadows, and the corresponding 
subset in $S_n$ will be denoted by $\bS(n)$. Then we have an induced partition 
$$\mathcal{A}=\mathcal{A}^0\cup\mathcal{A}^+\cup\mathcal{A}^-,$$ 
such that $\mathcal{A}^0\cup\mathcal{A}^+=\bS(n)$, $|\mathcal{A}^+|=|\mathcal{A}^-|$, for every matrix 
$A\in\mathcal{A}^0$, we have $[-A]=[A]$, and for any $A\in\mathcal{A}^+$, there is exactly one matrix 
$A^-\in\mathcal{A}^-$ with $[-A]=[A^-]$. \smallskip 

Of course, this set depends on the choice of $\mathcal{A}$, and whenever we use the term basic, we mean that we 
have a set constructed in the above way. In this paper, we give an algorithm which generates directly one particular 
set $\bS(n)$ of basic tame periodicity shadows.  \medskip 

For every other choice of the set $\mathcal{A}=\{\mathbb{A}_1',\dots,\mathbb{A}'_r\}$ of representatives of 
$G$-orbits in $S_n$, the resulting set $\bS'$ induces the same partition, up to exchanging matrices between 
the sets $\mathcal{A}^+$ and $\mathcal{A}^-$. More precisely, let 
$$\mathcal{A}^0=\{\mathbb{A}^0_1,\dots,\mathbb{A}^0_p\}, \ 
\mathcal{A}^+=\{\mathbb{A}^+_1,\dots,\mathbb{A}^+_q\}, \mbox { and } 
\mathcal{A}^-=\{\mathbb{A}^-_1,\dots,\mathbb{A}^-_q\}.$$ 
Then we can order matrices in $\bS'$ in such a way that $\bS'_k\sim \mathbb{A}^0_k$, for $k\in\{1,\dots,p\}$, 
$\bS_{p+k}'\sim \mathbb{A}^+_k$, for some $ 0 \leqslant q' \leqslant q$ and $k\in\{1,\dots,q'\}$, and $\bS_{p+k}' \sim -A^+_k$, 
if $q'+1 \leqslant k\leqslant q$. In particular, we counclude that $\bS'_k \sim \bS_k$, for any $k\in\{1,\dots,p+q'\}$ 
and $\bS'_k \sim -\bS_k$, for $k\in\{p+q'+1,\dots,p+q=s\}$. \medskip 

Summing up, we may work with a fixed set of basic tame periodicity shadows, since every other set basic shadows 
is obtained from $\bS(n)$ by taking permutations or opposite matrices. In particular, any such a set satisfies the 
following condition: {\it for any tame periodicity shadow $A\in S_n$ there exists a unique $i\in\{1,\dots,s\}$ such 
that $A \sim \bS_i$ or $A \sim -\bS_i$}. \medskip 

There is a related notion of an {\it essential shadow}, introduced to exclude some of the basic tame periodicity 
shadows, which cannot represent a Gabriel quiver of a tame symmetric algebra with periodic simples. We have a 
special shadow with double arrows, called a Markov shadow, which is the 
following matrix 
$$\bM=\left[\begin{array}{ccc} 0 & -2 & 2 \\ 2 & 0 & -2 \\ -2 & 2 & 0\end{array}\right].$$ 
Note that $\bM\in\bS(3)$ and the associated quiver $\bfQ_\bM$ is the Markov quiver (see Section \ref{sec:3}). \smallskip 

Now, a basic shadow $\bS=[a_{ij}]\in\bS(n)$ is called {\it essential}, if $\bS$ is the Markov shadow 
$\bS=\bM$, or $\bS\neq\bM$ satisfies conditions stated below. 
\begin{enumerate} 
\item[PS4)] Each row of $\bS$ does not contain both $2$ and $-2$.
\item[PS5)] For any $i,j,k$ such that $a_{ij}=2$ and $a_{jk}=1$, we have $a_{ki}>0$; for any $a_{ij}=-2$ and 
$a_{jk}=-1$, we have $a_{ki}<0$. 
\end{enumerate} \medskip 

Note that the condition PS4) is equivalent to say that the associated quiver $\bfQ_{\bS}$ has no vertex 
which is a source and a target of double arrows. If $\Lambda=KQ/I$ is a tame symmetric algebra with simple 
modules of period $4$ and $\bS_\La\neq\bM$, then its shadow $\bS_\La$ must satisfy PS4). This is a 
consequence of arguments presented in \cite[Lemma 5.2]{AGQT} which can be repeated in general. Alternatively, 
one can prove this using \cite[Lemmas 3.2 and 3.3]{Sko} and the fact that if $\La$ is tame, then given double 
arrows $\alpha,\ba:i\to j$ and an arrow $\beta:j\to k$ in $Q$, one of the paths $\alpha\beta$ or $\ba\beta$ 
must be involved in a minimal relation of $I$ \cite[see also Section 2]{Sko}. \smallskip 

Condition PS5) can be similarily rephrased in the language of quivers, namely, it is equivalent to say that the quiver 
$\bfQ_\bS$ has no subquiver of the form 
$$\xymatrix{i \ar@<+0.4ex>[r]^{\alpha} \ar@<-0.4ex>[r]_{\ba} & j \ar[r]^{\beta} & k} \qquad 
\mbox{or} \qquad \xymatrix{i \ar@<+0.4ex>@{<-}[r]^{\alpha} \ar@<-0.4ex>@{<-}[r]_{\ba} & j \ar@{<-}[r]^{\beta} & k}$$ 
such that there are no arrows $k\to i$ or $i\to k$, respectively. Now, it follows easily from \cite[Lemma 3.3]{Sko} 
that any such a subquiver gives two paths $\alpha\beta$ and $\ba\beta$ (or $\beta\alpha,\beta\ba$ in opposite case) 
not involved in minimal relations of $I$, and hence any algebra $\La=KQ/I$ whose shadow does not satisfy PS5) is a 
wild algebra. \smallskip 

Concluding, every algebra $\Lambda$ with $\bS_\Lambda$ not essential, is a wild algebra. So when we study tame 
algebras, we may restrict ourselves only to essential shadows. The lists of all essential shadows for $n\leqslant 6$ 
are presented in Sections \ref{sec:3}-\ref{sec:5}. Shadows which are not essential can be found in the Appendix. \medskip 

The algorithm is organised in the following way. We start with a related notion of a {\it shade}, by which we mean an 
integer matrix satisfying conditions PS1)-PS2) and T1)-T3). So the set of shades $\mathcal{S}_n$ contains set of all 
tame periodicity shadows $S_n$. The strategy is to first generate the set $\cS(n)$ of (basic) shades, and for any 
such a shade $A$, find all the solutions of the system of linear equations $AX=0$ and $X^T=X$. Then after checking if 
there exists a solution with $X\in\bM_n(\bN)$ and non-zero columns, we decide whether a given shade is a shadow, i.e. 
satisfies PS3). This gives the set $\bS(n)$ of all basic tame periodicity shadows. Final step is to run through the list 
$\bS(n)$, and check which of the shadows are essential. \medskip

Now, we are ready to present the algorithm. First, we introduce some preparatory notations. 
\smallskip

Let $\Sigma$ the set of permutations of the set $\{1,\dots,n\}$. For $A= (a_{ij})_{i,j\in\{1,\dots,n\}}\in S_n$ 
and $\sigma \in \Sigma$ we denote by $A_{\sigma}$ the matrix $(a_{\sigma(i)\sigma(j)})_{i,j\in\{1,\dots,n\}}$, which is 
exactly the matrix $P^TAP$, where $P=P_\sigma$ is the permutation matrix corresponding to $\sigma$.  
\smallskip

Further, we define the binary relation $\preceq$ on $S_n$ as follows. For any $A= (a_{ij})_{i,j\in\{1,\dots,n\}}\in S_n$ 
and $B= (b_{ij})_{i,j\in\{1,\dots,n\}}\in S_n$, we have $A \preceq B$ if and only if either $A=B$ or there are 
exist $i^*,j^*\in\{1,\dots,n\}$ satisfying the following conditions
\begin{enumerate}
  \item [(i)] $a_{ij} = b_{ij}$, for all $i\in\{1,\dots,i^*-1\}$ and $j\in\{1,\dots,n\}$,  
  \item [(ii)] $ a_{i^*j} = b_{i^*j}$, for all $j\in\{1,\dots,j^*-1\}$, 
  \item [(iii)] $a_{i^*j^*} < b_{i^*j^*}$. 
\end{enumerate} 
In the second case, we write $A\prec B$. \smallskip 

Clearly, the relation $\preceq$ is a partial oreder, and for any $A\neq B$ we have either $A\prec B$ or $B\prec A$, 
so in fact $S_n$ is totally ordered by the relation $\preceq$. Observe also that for each matrix $A \in S_n$ there 
exist exactly one matrix $B$, which is of the form $B=A_\sigma$, for some $\sigma\in\Sigma$, and  $B \preceq A_{\gamma}$, 
for any other $\gamma \in \Sigma$. Equivalently, one can treat $B$ as a matrix obtained from $A$ by permutations of rows 
and columns which is the smallest with respect to the order $\preceq$. We will denote such a matrix $B$ by $A_{\prec}$. 
Obviously, $A_\prec \preceq A$, by definition, because $A=A_{id}$ for the identity permutation $id$. Note also that the 
matrix $B$ is uniquely determined, but there may be many permutations $\sigma\in\Sigma$ realising $B=A_\sigma$. \smallskip 

For convenience, we define also matrix $A^*_\prec$ which is the largest element in the orbit $[A]$ with respect 
to the order $\preceq$. Clearly, we have $(-A)_\prec=-A_\prec^*$, $(A_\prec^*)_\prec=A_\prec$ and 
$(A_\prec)^*_\prec=A_\prec^*$. \smallskip 

We want to generate matrices satisfying conditions PS1), PS2), T1), T2), T3) recursively row by row. The main recursion 
takes as the input the size $n$ and the sequence $rows$ of the first $r$ rows of a matrix which is a candidate for a shadow, 
and recursively generate all shadows having these rows as the first $r$ rows. In each step the $r$ rows $R_1,\dots,R_r$ 
form a piece of potentially skew-symmetric matrix, by which we mean that the $r\times n$ matrix $R$ composed from these 
rows satisfies $R^T=-R$ at the admissible range of entries $(i,j)$, i.e. for $i,j\leqslant r$. In the main code we will 
use the following auxiliary functions. 

\begin{itemize}
  \item 
	\textsc{IsTPS}($row$) -- check if the given row satisfy the conditions T2), T3) and PS2); 
	it will follow from the construction that any such a row trivially satisfies T1); 
  \item 
	\textsc{ComposeRow}($rows$,$a_{r+2},\dots,a_n$) -- takes the sequence $rows$ of $r\in\{0,\dots,n-1\}$ rows 
	$R_1,\dots,R_r\in\bM_{1\times n}$, and additional $n-r-1$ new entries $a_k\in\{-2,-1,-0,1,2\}$, and returns 
	the the $(r+1)$-th row of the form $R_{r+1}=[\ \cdots 0 \ a_{r+2} \ \cdots \ a_n ]$ such that the new matrix 
	$\tilde{R}$ composed from rows $R_1,\dots,R_{r+1}$ also satisfies $\tilde{R}^T=-\tilde{R}$; note that in case $r=n-1$, 
	the last row $R_{r+1}=R_n$ is uniquely determined by the last column of matrix $R$, and the sequence of 
	additional entries is empty. 
  \item 
	\textsc{Matrix} -- construct the $n\times n$ matrix from the given sequence of rows. 
\end{itemize} \medskip  

Now, the main recursion is realised by the following procedure. \smallskip

	\begin{algorithmic}[1]
	\Procedure{SetRow}{$n, rows$}
	\State
	$r \gets \#rows$
	\If{$r = n-1$} 
	\State
	$lastRow \gets \textsc{ComposeRow}(rows)$  \Comment{the $n$-th row is given by the last column}
	\If {$\textsc{IsTPS}(lastRow)$}
	\State
	$rows \gets rows,lastRow$
	\State
	$M \gets \textsc{Matrix}(rows)$
	\If {$((2\!\not|\, n \vee \det(M) = 0) \wedge M_{\prec} = M \wedge M \preceq (-M)_{\prec}$}
	\State\Return $\{ M \}$ \Comment{$M$ is a new matrix}
	\Else
	\State\Return $\emptyset$
	\EndIf  
	\Else
	\State\Return $\emptyset$
	\EndIf   
	\Else
	\ForAll{$a_{r+2}, \dots, a_n \in \{-2,-1,0,1,2\}$}
	\Comment{coefficient which are not fixed by $rows$}
	\State
	$matrices \gets \emptyset$
	\State
	$nextRow \gets \textsc{ComposeRow}(rows,a_{r+2},...,a_n)$
	\Comment{compute the next row}
	\If {$\textsc{IsTPS}(nextrow)$}
	\State
	$newRows \gets rows,nextRow$
	\State
	$matrices \gets matrices \cup \textsc{SetRow}(n, newRows)$
	\EndIf  
	\EndFor
	\State\Return $matrices$
	\EndIf
	\EndProcedure
\end{algorithmic} \medskip 

Finally, applying the above procedure with initial data $(n,\emptyset)$, we get the following recursive algorithm  
generating the set $\cS(n)$ of all basic shades. \smallskip 

\begin{algorithmic}[H]
	\Procedure{MatricesSatisfyingTPS}{$n$}
	\State\Return $\textsc{SetRow}(n, \emptyset)$
	\EndProcedure
\end{algorithmic} \medskip

We note that the key line 8 of the code of \textsc{SetRow} allows us to remove matrices obtained by permuting rows 
and columns, or the opposite matrices $-M$. More precisely, if at some point we obtain a matrix $M$ with $M_\prec \neq M$, 
then $M_\prec \prec M$, so $M_\prec$ was counted to the set at some earlier step, and hence $M$ can be skipped, as 
it is obtained from $M_\prec $ by some permutation of rows and columns. Similarily, if $M\npreceq (-M)_{\prec}$, then (total 
order) $(-M)_\prec \prec M$, so $(-M)_\prec$ was counted earlier and $-M$ is its permutation. Observe also that the 
relation $\prec$ can be easily (and quickly) verified in an iterative manner. Moreover, the condition $M = M_{\prec}$ 
can be efficiently verified by a recursive procedure (we omit the details). \smallskip 

We also note that above procedure gives us the sequence of matrices ordered with respect to $\prec$. In particular, 
we can speed up the algorithm by stopping it when it generate a zero matrix. But the above code for simplicity omits 
this kind of optimization. \smallskip 

Let us also note that the presented algorithm can be easily parallelized because the iterations of the \textit{for} 
loop from lines 17-24 (apart from combining the obtained results into a list) are independent. This allows for better 
scaling of calculations. \medskip 

Now, it remains to show that the algorithm is correct, i.e. it indeed generates one particular set of basic 
shades. Denote by $\cS(n)$ the set generated by the procedure \textsc{MatricesSatisfyingTPS}. It is 
sufficient to prove that any shade belongs to $[S]$ or $[-S]$, for some $S\in\cS(n)$. Let $A$ be an 
arbitrary shade. \smallskip 

We have only one shade for $n=1$ or $2$, and this is the zero matrix, so we can assume that $n\geqslant 3$. 
Since we start with $rows=\emptyset$, the first iteration has $r=0<n-1$ and it initiates the loop generating 
all possible first rows $[a_1 \ a_2 \cdots \ a_n ]$. As a result, we can guarantee that the first row of $A$ 
is generated during this first step, and because $A$ is a shade, its first row satisfies T2), T3) and PS2), 
and hence, the function \textsc{IsTPS}($nextrow$) passes the first row to the next iteration. Similarily, the 
second row of $A$ is reached in the second iteration, and repeating this argument, we conclude that the matrix 
$M$ in the first if condition ($r=n-1$) is at some iteration equal to $A$. Because $A$ is a shade, the condition 
$2\!\not|n$ or $\det(A)=0$ always holds true. Finally, if $A_\prec =A$ and $A\preceq (-A)_\prec$, then $A$ 
is obviously in $\cS(n)$. If one of these conditions is not satisfied, then $A_\prec$ or $(-A)_\prec$ is 
less (and different) from $A$ in the order $\preceq$, and $A_\prec,(-A)_{\prec}\in S_n$, since these matrices 
are permutations of $A\in S_n$ or $-A\in S_n$. It means that both $A_\prec$ and $(-A)_{\prec}$ are shades, and 
hence they appear as $M$ in some other iteration. If $A_\prec \prec A$ and additionally $A_\prec \preceq (-A_\prec)_\prec$, 
then $A_\prec \in \cS(n)$, since we always have $(A_\prec)_\prec=A_\prec$. If $A_{\prec} \prec A$, but 
$(-A_\prec)_\prec \prec A_\prec $, then $(-A)_\prec \in\cS(n)$, because 
$$(-A)_\prec = -A^*_\prec = - (A_\prec)^*_\prec = (-A_\prec)_\prec \prec A_\prec = (A^*_\prec)_\prec=(-(-A)_\prec)_\prec.$$ 
Eventually, if $A_\prec=A$ but $(-A_\prec)_\prec \prec A_\prec$, then as above $(-A)_\prec\in\cS(n)$. Otherwise, 
we get $A_\prec=A\in\cS(n)$. $\hfill\Box$ \medskip 

As a result, we conclude that indeed the recursive procedure \textsc{MatricesSatisfyingTPS} generates a set 
of basic shades. Complete lists are presented in the Appendix. Due to their sizes, we present a full lists only 
for $n=3$ and $4$ (see Section \ref{sec:3}), whereas for $n=5$ and $6$, we restrict the lists only to essential 
shadows. Note that, each item on the lists in Sections \ref{sec:4} and \ref{sec:5} is a triple $(A,x,C)$, where 
$A$ is an essential periodicity shadow, $x$ is a generic vector of the nullspace $\mathcal{N}(A)$ (given in a 
parametric form), and $C$ is a generic matrix satisfying $AC=0$ and $C^T=C$. Moreover, since each $C$ is a symmetric 
matrix, we display only the coefficients over the main diagonal (including the main diagonal). In Section \ref{sec:3}, 
we skip $x$ and $C$ showing only shades (shadows); for the shape of $x$ and $C$, we refer interested reader to the 
Appendix, Section 2. \medskip 

Final step is to run through the list $\cS(n)$, and check whether a given shade $A$ is a tame periodicity 
shadow, i.e. it satisfies PS3). To do this, we exploit the form of solutions $C$ of linear system $AC=0$ and $C^T=C$, 
whose columns are vectors of the nullspace $\mathcal{N}(A)$. Indeed, to get the list of basic tame periodicity shadows 
$\bS(n)\subset\cS(n)$, we delete from $\cS(n)$ all the triples $(A,x,C)$, for which one of the following conditions 
hold 
\begin{itemize}
\item $x$ contains an entry equal to zero; 
\item $x$ contains two entries of the form $v$ and $-v$, where $v$ is a parameter; 
\item $x$ contains entries $\alpha_1,\dots,\alpha_r$ (not necesarily parameters) such that $a_1\alpha_1+\dots+a_r\alpha_r=0$, 
for some $a_1,\dots,a_r\in\bN$. 
\end{itemize} 

Note that each of these conditions implies that every solution $C$ either has a zero column or a negative entry. 
It is obvious in the first case, since then each column of $C$ has common entry equal to zero, so $C$ has a zero row, 
and hence a zero column, due to symmetricity. In the second case, we have either $v=-v=0$ in every column, and then 
we obtain a zero row (so a zero column), or there is one column for which $v>0$, but then the column contain 
a negative entry $-v$. In the last case, if $C$ has natural coefficients, we conclude that each of its columns has 
$\alpha_1=\dots=\alpha_r=0$, because $a_1,\dots,a_r\in\bN$, and then $C$ admits at least one zero row, hence a column. \medskip 

Finally, observe that after applying the above procedure, we obtain the set $\bS(n)$ of all basic tame periodicity shadows. 
Indeed, all the removed triples certainly not satisfy PS3), so they are not periodicity shadows. On the other hand, for 
any triple $(A,x,C)$ which left in $\bS(n)$, one can directly verify that there is at least one $C$ with coefficients in 
$\bN$ and non-zero columns. It requires to see that the condition $C\in\bM_n(\bN)$ is equivalent to a system of 
linear inequalities, which always have at least one solution with all $c_i\geqslant 0$ and all columns of $C$ non-zero. 
Verifying this can be done automatically via standard computational environments like Maple, but we skip this for 
simplicity. From the list of all basic tame periodicity shadows, one can easily get the list of essential shadows, 
by excluding matrices $A$ not satisfying PS4) and PS5). 

\bigskip

\section{Periodicity shadows of sizes $\leqslant$ 4}\label{sec:3} 

In this section we present the lists of tame periodicity shadows of small size $n\leqslant 4$. Since the cases 
$n=1$ or $2$ are trivial (there is only one tame periodicity shadow, which is the zero matrix), we will focus 
on the cases $n=3$ and $4$. \medskip 

Applying the algorithm from Section \ref{sec:2}, we obtain the list $\bS(3)$ containing $5$ matrices $\bS_1,\dots,\bS_5$, 
where the last one is zero, and the first four are the following matrices 
$$\left[\begin{array}{ccc} 0 & -2 & 2 \\ 2 & 0 & -2 \\ -2 & 2 & 0\end{array}\right]=\bM, \ 
\left[\begin{array}{ccc} 0 & -1 & 1 \\ 1 & 0 & -1 \\ -1 & 1 & 0\end{array}\right], \ 
\left[\begin{array}{ccc} 0 & -2 & 1 \\ 2 & 0 & -1 \\ -1 & 1 & 0\end{array}\right], \ 
\left[\begin{array}{ccc} 0 & -2 & 1 \\ 2 & 0 & -2 \\ -1 & 2 & 0\end{array}\right],$$ 
We will identify the above shadows with their corresponding quivers $\bbQ_1,\dots,\bbQ_5$, $\bbQ_i=\bfQ_{\bS_i}$, 
given as follows. 
$$\xymatrix@C=0.5cm{ &\circ \ar[rd]\ar@<+0.13cm>[rd] & \\ 
\circ \ar[ru]\ar@<+0.13cm>[ru] && \circ \ar[ll]\ar@<+0.13cm>[ll]} \ \ 
\xymatrix@C=0.5cm{ &\circ \ar[rd] & \\ \circ \ar[ru] && \circ \ar[ll]} \ \ 
\xymatrix@C=0.5cm{ &\circ \ar[rd] & \\ 
\circ \ar[ru] && \circ \ar[ll]\ar@<+0.13cm>[ll]} \ \ 
\xymatrix@C=0.5cm{ &\circ \ar[rd]\ar@<+0.13cm>[rd] & \\ 
\circ \ar[ru] && \circ \ar[ll]\ar@<+0.13cm>[ll]} \ \ 
\xymatrix@C=0.5cm{ &\circ  & \\ 
\circ && \circ }$$ 
Note that in this case shadows coincide with shades, i.e. $\bS(3)=\mathcal{S}(3)$. Moreover, the shadow $\bS_4$ 
is not essential, since it does not satisfy condition PS4). Hence we have $4$ essential shadows: $3$ non-trivial 
$\bbQ_1,\bbQ_2,\bbQ_3$ and $\bbQ_5=\emptyset$. \medskip

After running the algorithm for $n=4$, we conclude that the list $\bS(4)=\mathcal{S}(4)$ consists of the following 
$12$ matrices  $\bS_1,\dots,\bS_{12}$: 
$$\begin{bmatrix}
     0 & -2 &  0 &  1 \\
     2 &  0 &  0 & -1 \\
     0 &  0 &  0 &  0 \\
    -1 &  1 &  0 &  0 \\
  \end{bmatrix}, \ \begin{bmatrix}
     0 & -1 &  0 &  1 \\
     1 &  0 & -1 &  0 \\
     0 &  1 &  0 & -1 \\
    -1 &  0 &  1 &  0 \\
  \end{bmatrix}, \ \begin{bmatrix}
     0 & -1 &  0 &  1 \\
     1 &  0 &  0 & -1 \\
     0 &  0 &  0 &  0 \\
    -1 &  1 &  0 &  0 \\
  \end{bmatrix}, \  \begin{bmatrix}
     0 & -2 &  1 &  1 \\
     2 &  0 & -1 & -1 \\
    -1 &  1 &  0 &  0 \\
    -1 &  1 &  0 &  0 \\
  \end{bmatrix},$$ 
  
$$\begin{bmatrix}
     0 & -1 & -1 &  1 \\
     1 &  0 & -1 &  0 \\
     1 &  1 &  0 & -1 \\
    -1 &  0 &  1 &  0 \\
  \end{bmatrix}, \ \begin{bmatrix}
     0 & -1 & -1 &  1 \\
     1 &  0 &  0 & -1 \\
     1 &  0 &  0 & -1 \\
    -1 &  1 &  1 &  0 \\
  \end{bmatrix}, \ \begin{bmatrix}
     0 &  0 &  0 &  0 \\
     0 &  0 &  0 &  0 \\
     0 &  0 &  0 &  0 \\
     0 &  0 &  0 &  0 \\
  \end{bmatrix}, \ \begin{bmatrix}
     0 & -2 &  0 &  1 \\
     2 &  0 & -2 &  0 \\
     0 &  2 &  0 & -1 \\
    -1 &  0 &  1 &  0 \\
  \end{bmatrix}$$ 
  
$$\begin{bmatrix}
     0 & -2 &  0 &  1 \\
     2 &  0 &  0 & -2 \\
     0 &  0 &  0 &  0 \\
    -1 &  2 &  0 &  0 \\
  \end{bmatrix}, \ \begin{bmatrix}
     0 & -2 &  0 &  2 \\
     2 &  0 & -2 &  0 \\
     0 &  2 &  0 & -2 \\
    -2 &  0 &  2 &  0 \\
  \end{bmatrix}, \ \begin{bmatrix}
     0 & -2 &  0 &  2 \\
     2 &  0 &  0 & -2 \\
     0 &  0 &  0 &  0 \\
    -2 &  2 &  0 &  0 \\
  \end{bmatrix}, \ \begin{bmatrix}
     0 & -2 &  0 &  2 \\
     2 &  0 & -1 & -1 \\
     0 &  1 &  0 & -1 \\
    -2 &  1 &  1 &  0 \\
  \end{bmatrix}.$$

The corresponding quivers $\bbQ_1,\dots,\bbQ_{12}$, $\bbQ_i=\bfQ_{\bS_i}$, are given as follows.  

$$\xymatrix@C=0.5cm@R=0.5cm{\circ \ar[r] & \circ \ar[ld] \\ 
\circ \ar[u]\ar@<+0.1cm>[u]& \circ } \ 
\xymatrix@C=0.5cm@R=0.5cm{\circ \ar[r] & \circ \ar[d] \\ 
\circ \ar[u] & \circ \ar[l] } 
\xymatrix@C=0.5cm@R=0.5cm{\circ \ar[r] & \circ \ar[ld] \\ 
\circ \ar[u] & \circ } \ 
\xymatrix@C=0.5cm@R=0.5cm{\circ \ar[r]\ar[rd] & \circ \ar[ld] \\ 
\circ \ar[u]\ar@<+0.1cm>[u] & \circ \ar[l] } \ 
\xymatrix@C=0.5cm@R=0.5cm{\circ \ar[r] & \circ \ar[d] \\ 
\circ \ar[u] & \circ \ar[l]\ar[lu] } \ 
\xymatrix@C=0.5cm@R=0.5cm{\circ \ar[r] & \circ \ar[d]\ar[ld] \\ 
\circ \ar[u] & \circ \ar[lu] }$$ 

$$\xymatrix@C=0.5cm@R=0.5cm{\circ & \circ \\ 
\circ & \circ}
\xymatrix@C=0.5cm@R=0.5cm{\circ \ar[r] & \circ \ar[d] \\ 
\circ \ar[u]\ar@<+0.1cm>[u]& \circ \ar[l]\ar@<+0.1cm>[l]} \ 
\xymatrix@C=0.5cm@R=0.5cm{\circ \ar[r] & \circ \ar[ld]\ar@<+0.1cm>[ld] \\ 
\circ \ar[u]\ar@<+0.1cm>[u]& \circ } \ \xymatrix@C=0.5cm@R=0.5cm{\circ \ar[r]\ar@<+0.1cm>[r] & \circ \ar[d]\ar@<+0.1cm>[d] \\ 
\circ \ar[u]\ar@<+0.1cm>[u]& \circ \ar[l]\ar@<+0.1cm>[l]} \ 
\xymatrix@C=0.5cm@R=0.5cm{\circ \ar[r]\ar@<+0.1cm>[r] & \circ \ar[ld]\ar@<+0.1cm>[ld] \\ 
\circ \ar[u]\ar@<+0.1cm>[u]& \circ} \ 
\xymatrix@C=0.5cm@R=0.5cm{\circ \ar[r]\ar@<+0.1cm>[r] & \circ \ar[d]\ar[ld] \\ 
\circ \ar[u]\ar@<+0.1cm>[u]& \circ \ar[l]}$$ \smallskip 
 
It is easy to see that shadows $\bS_i$, for $i\geqslant 8$ do not satisfy PS4). Additionally, the shadows 
$\bS_8$ and $\bS_{10}$ do not satisfy PS5). Thus we get $7$ essential shadows, which correspond to seven 
essential quivers $\bbQ_1,\dots,\bbQ_{7}$ (including the empty one). In particular, we conclude that 
for every algebra $\La=KQ/I$ of generalized quaternion type with $Q=Q_\La$ having four vertices, its Gabriel 
quiver $Q$ is obtained (modulo loops) from the quivers $\bbQ_1,\dots,\bbQ_7$ by adding $2$-cycles in a very 
restrictive way (as described in the main result \cite{Sko}). Due to \cite[Corollary 5.3]{Sko}, we may exclude 
the empty shadow $\bbQ_7$, since it has too many vertices. \medskip   

The following table summarizes the numbers of shades and shadows (including essential ones) in the cases 
$n=3,4,5$ and $6$. \medskip 

\begin{tabular}{|p{2cm}|p{4cm}|p{4cm}|p{4cm}|} 
\hline 
& number of shades & number of shadows & essential shadows \\ 
\hline 
n=3 & 5 & 5 & 4 \\ 
n=4 & 12 & 12 & 7 \\ 
n=5 & 138 & 65 & 26 \\ 
n=6 & 1260 & 516 & 223 \\ 
\hline 
\end{tabular}

\bigskip

\section{Periodicity shadows of size n=5}\label{sec:4} 

In this section we present the list of all essential shadows of size $n=5$. In this case, we have exactly $26$ 
essential shadows among all $65$ (basic) tame periodicity shadows in $\bS(5)$. The remaining non-essential shadows and 
the complete list of all $|\cS(5)|=138$ (basic) shades, including those which are not tame periodicity shadows, 
can be found in the Appendix. \medskip


\parbox{.06\linewidth}{\stepcounter{mac}(\themac)}%
\resizebox{.43\linewidth}{!}{%
$
  \begin{bmatrix}
     0 & -2 &  0 &  0 &  1 \\
     2 &  0 &  0 &  0 & -1 \\
     0 &  0 &  0 & -2 &  1 \\
     0 &  0 &  2 &  0 & -1 \\
    -1 &  1 & -1 &  1 &  0 \\
  \end{bmatrix}
\begin{bmatrix} v \\ v \\ v \\ v \\  2  v  \end{bmatrix}
  \begin{bmatrix}
    c & c & c & c &  2  c  \\
    \cdot & c & c & c &  2  c  \\
    \cdot & \cdot & c & c &  2  c  \\
    \cdot & \cdot & \cdot & c &  2  c  \\
    \cdot & \cdot & \cdot & \cdot &  4  c  \\
  \end{bmatrix}
$
}
\makebox[.01\linewidth][r]{\,}%
\parbox{.06\linewidth}{\stepcounter{mac}(\themac)}%
\resizebox{.43\linewidth}{!}{%
$
  \begin{bmatrix}
     0 & -2 &  0 &  0 &  1 \\
     2 &  0 &  0 &  0 & -1 \\
     0 &  0 &  0 & -1 &  1 \\
     0 &  0 &  1 &  0 & -1 \\
    -1 &  1 & -1 &  1 &  0 \\
  \end{bmatrix}
\begin{bmatrix} v \\ v \\  2  v  \\  2  v  \\  2  v  \end{bmatrix}
  \begin{bmatrix}
    c & c &  2  c  &  2  c  &  2  c  \\
    \cdot & c &  2  c  &  2  c  &  2  c  \\
    \cdot & \cdot &  4  c  &  4  c  &  4  c  \\
    \cdot & \cdot & \cdot &  4  c  &  4  c  \\
    \cdot & \cdot & \cdot & \cdot &  4  c  \\
  \end{bmatrix}
$
}

\parbox{.06\linewidth}{\stepcounter{mac}(\themac)}%
\resizebox{.43\linewidth}{!}{%
$
  \begin{bmatrix}
     0 & -2 &  0 &  0 &  1 \\
     2 &  0 &  0 &  0 & -1 \\
     0 &  0 &  0 &  0 &  0 \\
     0 &  0 &  0 &  0 &  0 \\
    -1 &  1 &  0 &  0 &  0 \\
  \end{bmatrix}
\begin{bmatrix} v_{1} \\ v_{1} \\ v_{2} \\ v_{3} \\  2  v_{1}  \end{bmatrix}
  \begin{bmatrix}
    c_{1} & c_{1} & c_{2} & c_{3} &  2  c_{1}  \\
    \cdot & c_{1} & c_{2} & c_{3} &  2  c_{1}  \\
    \cdot & \cdot & c_{4} & c_{5} &  2  c_{2}  \\
    \cdot & \cdot & \cdot & c_{6} &  2  c_{3}  \\
    \cdot & \cdot & \cdot & \cdot &  4  c_{1}  \\
  \end{bmatrix}
$
}
\makebox[.01\linewidth][r]{\,}%
\parbox{.06\linewidth}{\stepcounter{mac}(\themac)}%
\resizebox{.43\linewidth}{!}{%
$
  \begin{bmatrix}
     0 & -2 &  0 &  1 &  1 \\
     2 &  0 &  0 & -1 & -1 \\
     0 &  0 &  0 & -1 &  1 \\
    -1 &  1 &  1 &  0 & -1 \\
    -1 &  1 & -1 &  1 &  0 \\
  \end{bmatrix}
\begin{bmatrix} v \\ v \\ v \\ v \\ v \end{bmatrix}
  \begin{bmatrix}
    c & c & c & c & c \\
    \cdot & c & c & c & c \\
    \cdot & \cdot & c & c & c \\
    \cdot & \cdot & \cdot & c & c \\
    \cdot & \cdot & \cdot & \cdot & c \\
  \end{bmatrix}
$
}

\parbox{.06\linewidth}{\stepcounter{mac}(\themac)}%
\resizebox{.43\linewidth}{!}{%
$
  \begin{bmatrix}
     0 & -2 &  1 &  1 &  1 \\
     2 &  0 & -1 & -1 & -1 \\
    -1 &  1 &  0 & -1 &  1 \\
    -1 &  1 &  1 &  0 & -1 \\
    -1 &  1 & -1 &  1 &  0 \\
  \end{bmatrix}
\begin{bmatrix} \frac{ 3 }{ 2 } v  \\ \frac{ 3 }{ 2 } v  \\ v \\ v \\ v \end{bmatrix}
  \begin{bmatrix}
    \frac{ 9 }{ 4 } c  & \frac{ 9 }{ 4 } c  & \frac{ 3 }{ 2 } c  & \frac{ 3 }{ 2 } c  & \frac{ 3 }{ 2 } c  \\
    \cdot & \frac{ 9 }{ 4 } c  & \frac{ 3 }{ 2 } c  & \frac{ 3 }{ 2 } c  & \frac{ 3 }{ 2 } c  \\
    \cdot & \cdot & c & c & c \\
    \cdot & \cdot & \cdot & c & c \\
    \cdot & \cdot & \cdot & \cdot & c \\
  \end{bmatrix}
$
}
\makebox[.01\linewidth][r]{\,}%
\parbox{.06\linewidth}{\stepcounter{mac}(\themac)}%
\resizebox{.43\linewidth}{!}{%
$
  \begin{bmatrix}
     0 & -1 & -1 & -1 &  1 \\
     1 &  0 & -1 &  1 & -1 \\
     1 &  1 &  0 & -1 & -1 \\
     1 & -1 &  1 &  0 & -1 \\
    -1 &  1 &  1 &  1 &  0 \\
  \end{bmatrix}
\begin{bmatrix}  3  v  \\ v \\ v \\ v \\  3  v  \end{bmatrix}
  \begin{bmatrix}
     9  c  &  3  c  &  3  c  &  3  c  &  9  c  \\
    \cdot & c & c & c &  3  c  \\
    \cdot & \cdot & c & c &  3  c  \\
    \cdot & \cdot & \cdot & c &  3  c  \\
    \cdot & \cdot & \cdot & \cdot &  9  c  \\
  \end{bmatrix}
$
}

\parbox{.06\linewidth}{\stepcounter{mac}(\themac)}%
\resizebox{.43\linewidth}{!}{%
$
  \begin{bmatrix}
     0 & -1 & -1 &  0 &  1 \\
     1 &  0 & -1 &  1 & -1 \\
     1 &  1 &  0 & -1 & -1 \\
     0 & -1 &  1 &  0 &  0 \\
    -1 &  1 &  1 &  0 &  0 \\
  \end{bmatrix}
\begin{bmatrix}  2  v  \\ v \\ v \\ v \\  2  v  \end{bmatrix}
  \begin{bmatrix}
     4  c  &  2  c  &  2  c  &  2  c  &  4  c  \\
    \cdot & c & c & c &  2  c  \\
    \cdot & \cdot & c & c &  2  c  \\
    \cdot & \cdot & \cdot & c &  2  c  \\
    \cdot & \cdot & \cdot & \cdot &  4  c  \\
  \end{bmatrix}
$
}
\makebox[.01\linewidth][r]{\,}%
\parbox{.06\linewidth}{\stepcounter{mac}(\themac)}%
\resizebox{.43\linewidth}{!}{%
$
  \begin{bmatrix}
     0 & -1 & -1 &  0 &  1 \\
     1 &  0 & -1 &  1 & -1 \\
     1 &  1 &  0 & -1 &  0 \\
     0 & -1 &  1 &  0 &  0 \\
    -1 &  1 &  0 &  0 &  0 \\
  \end{bmatrix}
\begin{bmatrix} v \\ v \\ v \\  2  v  \\  2  v  \end{bmatrix}
  \begin{bmatrix}
    c & c & c &  2  c  &  2  c  \\
    \cdot & c & c &  2  c  &  2  c  \\
    \cdot & \cdot & c &  2  c  &  2  c  \\
    \cdot & \cdot & \cdot &  4  c  &  4  c  \\
    \cdot & \cdot & \cdot & \cdot &  4  c  \\
  \end{bmatrix}
$
}

\parbox{.06\linewidth}{\stepcounter{mac}(\themac)}%
\resizebox{.43\linewidth}{!}{%
$
  \begin{bmatrix}
     0 & -1 & -1 &  0 &  1 \\
     1 &  0 &  0 & -1 &  0 \\
     1 &  0 &  0 &  1 & -1 \\
     0 &  1 & -1 &  0 &  0 \\
    -1 &  0 &  1 &  0 &  0 \\
  \end{bmatrix}
\begin{bmatrix} v \\ v \\ v \\ v \\  2  v  \end{bmatrix}
  \begin{bmatrix}
    c & c & c & c &  2  c  \\
    \cdot & c & c & c &  2  c  \\
    \cdot & \cdot & c & c &  2  c  \\
    \cdot & \cdot & \cdot & c &  2  c  \\
    \cdot & \cdot & \cdot & \cdot &  4  c  \\
  \end{bmatrix}
$
}
\makebox[.01\linewidth][r]{\,}%
\parbox{.06\linewidth}{\stepcounter{mac}(\themac)}%
\resizebox{.43\linewidth}{!}{%
$
  \begin{bmatrix}
     0 & -1 & -1 &  1 &  1 \\
     1 &  0 & -1 & -1 &  1 \\
     1 &  1 &  0 & -1 & -1 \\
    -1 &  1 &  1 &  0 & -1 \\
    -1 & -1 &  1 &  1 &  0 \\
  \end{bmatrix}
\begin{bmatrix} v \\ v \\ v \\ v \\ v \end{bmatrix}
  \begin{bmatrix}
    c & c & c & c & c \\
    \cdot & c & c & c & c \\
    \cdot & \cdot & c & c & c \\
    \cdot & \cdot & \cdot & c & c \\
    \cdot & \cdot & \cdot & \cdot & c \\
  \end{bmatrix}
$
}

\parbox{.06\linewidth}{\stepcounter{mac}(\themac)}%
\resizebox{.43\linewidth}{!}{%
$
  \begin{bmatrix}
     0 & -1 & -1 &  1 &  1 \\
     1 &  0 &  0 & -1 &  0 \\
     1 &  0 &  0 &  0 & -1 \\
    -1 &  1 &  0 &  0 &  0 \\
    -1 &  0 &  1 &  0 &  0 \\
  \end{bmatrix}
\begin{bmatrix} v \\ v \\ v \\ v \\ v \end{bmatrix}
  \begin{bmatrix}
    c & c & c & c & c \\
    \cdot & c & c & c & c \\
    \cdot & \cdot & c & c & c \\
    \cdot & \cdot & \cdot & c & c \\
    \cdot & \cdot & \cdot & \cdot & c \\
  \end{bmatrix}
$
}
\makebox[.01\linewidth][r]{\,}%
\parbox{.06\linewidth}{\stepcounter{mac}(\themac)}%
\resizebox{.43\linewidth}{!}{%
$
  \begin{bmatrix}
     0 & -1 &  0 &  0 &  1 \\
     1 &  0 & -1 &  0 &  0 \\
     0 &  1 &  0 & -1 &  0 \\
     0 &  0 &  1 &  0 & -1 \\
    -1 &  0 &  0 &  1 &  0 \\
  \end{bmatrix}
\begin{bmatrix} v \\ v \\ v \\ v \\ v \end{bmatrix}
  \begin{bmatrix}
    c & c & c & c & c \\
    \cdot & c & c & c & c \\
    \cdot & \cdot & c & c & c \\
    \cdot & \cdot & \cdot & c & c \\
    \cdot & \cdot & \cdot & \cdot & c \\
  \end{bmatrix}
$
}

\parbox{.06\linewidth}{\stepcounter{mac}(\themac)}%
\resizebox{.43\linewidth}{!}{%
$
  \begin{bmatrix}
     0 & -1 &  0 &  0 &  1 \\
     1 &  0 & -1 &  0 &  0 \\
     0 &  1 &  0 &  0 & -1 \\
     0 &  0 &  0 &  0 &  0 \\
    -1 &  0 &  1 &  0 &  0 \\
  \end{bmatrix}
\begin{bmatrix} v_{1} \\ v_{3} \\ v_{1} \\ v_{2} \\ v_{3} \end{bmatrix}
  \begin{bmatrix}
    c_{1} & c_{3} & c_{1} & c_{2} & c_{3} \\
    \cdot & c_{6} & c_{3} & c_{5} & c_{6} \\
    \cdot & \cdot & c_{1} & c_{2} & c_{3} \\
    \cdot & \cdot & \cdot & c_{4} & c_{5} \\
    \cdot & \cdot & \cdot & \cdot & c_{6} \\
  \end{bmatrix}
$
}
\makebox[.01\linewidth][r]{\,}%
\parbox{.06\linewidth}{\stepcounter{mac}(\themac)}%
\resizebox{.43\linewidth}{!}{%
$
  \begin{bmatrix}
     0 & -1 &  0 &  0 &  1 \\
     1 &  0 &  0 &  0 & -1 \\
     0 &  0 &  0 &  0 &  0 \\
     0 &  0 &  0 &  0 &  0 \\
    -1 &  1 &  0 &  0 &  0 \\
  \end{bmatrix}
\begin{bmatrix} v_{3} \\ v_{3} \\ v_{1} \\ v_{2} \\ v_{3} \end{bmatrix}
  \begin{bmatrix}
    c_{6} & c_{6} & c_{3} & c_{5} & c_{6} \\
    \cdot & c_{6} & c_{3} & c_{5} & c_{6} \\
    \cdot & \cdot & c_{1} & c_{2} & c_{3} \\
    \cdot & \cdot & \cdot & c_{4} & c_{5} \\
    \cdot & \cdot & \cdot & \cdot & c_{6} \\
  \end{bmatrix}
$
}

\parbox{.06\linewidth}{\stepcounter{mac}(\themac)}%
\resizebox{.43\linewidth}{!}{%
$
  \begin{bmatrix}
     0 &  0 &  0 &  0 &  0 \\
     0 &  0 &  0 &  0 &  0 \\
     0 &  0 &  0 &  0 &  0 \\
     0 &  0 &  0 &  0 &  0 \\
     0 &  0 &  0 &  0 &  0 \\
  \end{bmatrix}
\begin{bmatrix} v_{1} \\ v_{2} \\ v_{3} \\ v_{4} \\ v_{5} \end{bmatrix}
  \begin{bmatrix}
    c_{1} & c_{2} & c_{3} & c_{4} & c_{5} \\
    \cdot & c_{6} & c_{7} & c_{8} & c_{9} \\
    \cdot & \cdot & c_{10} & c_{11} & c_{12} \\
    \cdot & \cdot & \cdot & c_{13} & c_{14} \\
    \cdot & \cdot & \cdot & \cdot & c_{15} \\
  \end{bmatrix}
$
}
\makebox[.01\linewidth][r]{\,}%

\parbox{.06\linewidth}{\stepcounter{mac}(\themac)}%
\resizebox{.94\linewidth}{!}{%
$
  \begin{bmatrix}
     0 & -2 &  0 &  1 &  1 \\
     2 &  0 &  0 & -1 & -1 \\
     0 &  0 &  0 &  0 &  0 \\
    -1 &  1 &  0 &  0 &  0 \\
    -1 &  1 &  0 &  0 &  0 \\
  \end{bmatrix}
\begin{bmatrix} v_{1} \\ v_{1} \\ v_{2} \\  - v_{3}   +  2  v_{1}   \\ v_{3} \end{bmatrix}
  \begin{bmatrix}
    c_{1} & c_{1} & c_{2} &  - c_{3}   +  2  c_{1}   & c_{3} \\
    \cdot & c_{1} & c_{2} &  - c_{3}   +  2  c_{1}   & c_{3} \\
    \cdot & \cdot & c_{4} &  - c_{5}   +  2  c_{2}   & c_{5} \\
    \cdot & \cdot & \cdot & c_{6}  -4  c_{3}   +  4  c_{1}   &  - c_{6}   +  2  c_{3}   \\
    \cdot & \cdot & \cdot & \cdot & c_{6} \\
  \end{bmatrix}
$
}

\parbox{.06\linewidth}{\stepcounter{mac}(\themac)}%
\resizebox{.94\linewidth}{!}{%
$
  \begin{bmatrix}
     0 & -2 &  1 &  1 &  1 \\
     2 &  0 & -1 & -1 & -1 \\
    -1 &  1 &  0 &  0 &  0 \\
    -1 &  1 &  0 &  0 &  0 \\
    -1 &  1 &  0 &  0 &  0 \\
  \end{bmatrix}
\begin{bmatrix} v_{1} \\ v_{1} \\  - v_{3}   - v_{2}   +  2  v_{1}   \\ v_{2} \\ v_{3} \end{bmatrix}
  \begin{bmatrix}
    c_{1} & c_{1} &  - c_{3}   - c_{2}   +  2  c_{1}   & c_{2} & c_{3} \\
    \cdot & c_{1} &  - c_{3}   - c_{2}   +  2  c_{1}   & c_{2} & c_{3} \\
    \cdot & \cdot & c_{6}  +  2  c_{5}   -4  c_{3}   + c_{4}  -4  c_{2}   +  4  c_{1}   &  - c_{5}   - c_{4}   +  2  c_{2}   &  - c_{6}   - c_{5}   +  2  c_{3}   \\
    \cdot & \cdot & \cdot & c_{4} & c_{5} \\
    \cdot & \cdot & \cdot & \cdot & c_{6} \\
  \end{bmatrix}
$
}

\parbox{.06\linewidth}{\stepcounter{mac}(\themac)}%
\resizebox{.94\linewidth}{!}{%
$
  \begin{bmatrix}
     0 & -1 & -1 & -1 &  1 \\
     1 &  0 & -1 & -1 &  0 \\
     1 &  1 &  0 &  0 & -1 \\
     1 &  1 &  0 &  0 & -1 \\
    -1 &  0 &  1 &  1 &  0 \\
  \end{bmatrix}
\begin{bmatrix} v_{2}  + v_{1}  \\ v_{3}  - v_{2}   - v_{1}   \\ v_{1} \\ v_{2} \\ v_{3} \end{bmatrix}
  \begin{bmatrix}
    c_{4}  +  2  c_{2}   + c_{1}  & c_{5}  - c_{4}   -2  c_{2}   + c_{3}  - c_{1}   & c_{2}  + c_{1}  & c_{4}  + c_{2}  & c_{5}  + c_{3}  \\
    \cdot & c_{6}  -2  c_{5}   -2  c_{3}   + c_{4}  +  2  c_{2}   + c_{1}  & c_{3}  - c_{2}   - c_{1}   & c_{5}  - c_{4}   - c_{2}   & c_{6}  - c_{5}   - c_{3}   \\
    \cdot & \cdot & c_{1} & c_{2} & c_{3} \\
    \cdot & \cdot & \cdot & c_{4} & c_{5} \\
    \cdot & \cdot & \cdot & \cdot & c_{6} \\
  \end{bmatrix}
$
}

\parbox{.06\linewidth}{\stepcounter{mac}(\themac)}%
\resizebox{.94\linewidth}{!}{%
$
  \begin{bmatrix}
     0 & -1 & -1 & -1 &  1 \\
     1 &  0 & -1 &  0 &  0 \\
     1 &  1 &  0 &  1 & -1 \\
     1 &  0 & -1 &  0 &  0 \\
    -1 &  0 &  1 &  0 &  0 \\
  \end{bmatrix}
\begin{bmatrix} v_{1} \\ v_{3}  - v_{2}   - v_{1}   \\ v_{1} \\ v_{2} \\ v_{3} \end{bmatrix}
  \begin{bmatrix}
    c_{1} & c_{3}  - c_{2}   - c_{1}   & c_{1} & c_{2} & c_{3} \\
    \cdot & c_{6}  -2  c_{5}   -2  c_{3}   + c_{4}  +  2  c_{2}   + c_{1}  & c_{3}  - c_{2}   - c_{1}   & c_{5}  - c_{4}   - c_{2}   & c_{6}  - c_{5}   - c_{3}   \\
    \cdot & \cdot & c_{1} & c_{2} & c_{3} \\
    \cdot & \cdot & \cdot & c_{4} & c_{5} \\
    \cdot & \cdot & \cdot & \cdot & c_{6} \\
  \end{bmatrix}
$
}

\parbox{.06\linewidth}{\stepcounter{mac}(\themac)}%
\resizebox{.94\linewidth}{!}{%
$
  \begin{bmatrix}
     0 & -1 & -1 & -1 &  1 \\
     1 &  0 &  0 &  0 & -1 \\
     1 &  0 &  0 &  0 & -1 \\
     1 &  0 &  0 &  0 & -1 \\
    -1 &  1 &  1 &  1 &  0 \\
  \end{bmatrix}
\begin{bmatrix} v_{3} \\ v_{3}  - v_{2}   - v_{1}   \\ v_{1} \\ v_{2} \\ v_{3} \end{bmatrix}
  \begin{bmatrix}
    c_{6} & c_{6}  - c_{5}   - c_{3}   & c_{3} & c_{5} & c_{6} \\
    \cdot & c_{6}  -2  c_{5}   -2  c_{3}   + c_{4}  +  2  c_{2}   + c_{1}  & c_{3}  - c_{2}   - c_{1}   & c_{5}  - c_{4}   - c_{2}   & c_{6}  - c_{5}   - c_{3}   \\
    \cdot & \cdot & c_{1} & c_{2} & c_{3} \\
    \cdot & \cdot & \cdot & c_{4} & c_{5} \\
    \cdot & \cdot & \cdot & \cdot & c_{6} \\
  \end{bmatrix}
$
}

\parbox{.06\linewidth}{\stepcounter{mac}(\themac)}%
\resizebox{.94\linewidth}{!}{%
$
  \begin{bmatrix}
     0 & -1 & -1 &  0 &  1 \\
     1 &  0 & -1 & -1 &  0 \\
     1 &  1 &  0 & -1 & -1 \\
     0 &  1 &  1 &  0 & -1 \\
    -1 &  0 &  1 &  1 &  0 \\
  \end{bmatrix}
\begin{bmatrix} v_{2}  + v_{1}  \\ v_{3}  - v_{1}   \\ v_{1} \\ v_{2} \\ v_{3} \end{bmatrix}
  \begin{bmatrix}
    c_{4}  +  2  c_{2}   + c_{1}  & c_{5}  + c_{3}  - c_{2}   - c_{1}   & c_{2}  + c_{1}  & c_{4}  + c_{2}  & c_{5}  + c_{3}  \\
    \cdot & c_{6}  -2  c_{3}   + c_{1}  & c_{3}  - c_{1}   & c_{5}  - c_{2}   & c_{6}  - c_{3}   \\
    \cdot & \cdot & c_{1} & c_{2} & c_{3} \\
    \cdot & \cdot & \cdot & c_{4} & c_{5} \\
    \cdot & \cdot & \cdot & \cdot & c_{6} \\
  \end{bmatrix}
$
}

\parbox{.06\linewidth}{\stepcounter{mac}(\themac)}%
\resizebox{.94\linewidth}{!}{%
$
  \begin{bmatrix}
     0 & -1 & -1 &  0 &  1 \\
     1 &  0 & -1 &  0 &  0 \\
     1 &  1 &  0 &  0 & -1 \\
     0 &  0 &  0 &  0 &  0 \\
    -1 &  0 &  1 &  0 &  0 \\
  \end{bmatrix}
\begin{bmatrix} v_{1} \\ v_{3}  - v_{1}   \\ v_{1} \\ v_{2} \\ v_{3} \end{bmatrix}
  \begin{bmatrix}
    c_{1} & c_{3}  - c_{1}   & c_{1} & c_{2} & c_{3} \\
    \cdot & c_{6}  -2  c_{3}   + c_{1}  & c_{3}  - c_{1}   & c_{5}  - c_{2}   & c_{6}  - c_{3}   \\
    \cdot & \cdot & c_{1} & c_{2} & c_{3} \\
    \cdot & \cdot & \cdot & c_{4} & c_{5} \\
    \cdot & \cdot & \cdot & \cdot & c_{6} \\
  \end{bmatrix}
$
}

\parbox{.06\linewidth}{\stepcounter{mac}(\themac)}%
\resizebox{.94\linewidth}{!}{%
$
  \begin{bmatrix}
     0 & -1 & -1 &  0 &  1 \\
     1 &  0 &  0 & -1 &  0 \\
     1 &  0 &  0 & -1 &  0 \\
     0 &  1 &  1 &  0 & -1 \\
    -1 &  0 &  0 &  1 &  0 \\
  \end{bmatrix}
\begin{bmatrix} v_{2} \\ v_{3}  - v_{1}   \\ v_{1} \\ v_{2} \\ v_{3} \end{bmatrix}
  \begin{bmatrix}
    c_{4} & c_{5}  - c_{2}   & c_{2} & c_{4} & c_{5} \\
    \cdot & c_{6}  -2  c_{3}   + c_{1}  & c_{3}  - c_{1}   & c_{5}  - c_{2}   & c_{6}  - c_{3}   \\
    \cdot & \cdot & c_{1} & c_{2} & c_{3} \\
    \cdot & \cdot & \cdot & c_{4} & c_{5} \\
    \cdot & \cdot & \cdot & \cdot & c_{6} \\
  \end{bmatrix}
$
}

\parbox{.06\linewidth}{\stepcounter{mac}(\themac)}%
\resizebox{.94\linewidth}{!}{%
$
  \begin{bmatrix}
     0 & -1 & -1 &  0 &  1 \\
     1 &  0 &  0 &  0 & -1 \\
     1 &  0 &  0 &  0 & -1 \\
     0 &  0 &  0 &  0 &  0 \\
    -1 &  1 &  1 &  0 &  0 \\
  \end{bmatrix}
\begin{bmatrix} v_{3} \\ v_{3}  - v_{1}   \\ v_{1} \\ v_{2} \\ v_{3} \end{bmatrix}
  \begin{bmatrix}
    c_{6} & c_{6}  - c_{3}   & c_{3} & c_{5} & c_{6} \\
    \cdot & c_{6}  -2  c_{3}   + c_{1}  & c_{3}  - c_{1}   & c_{5}  - c_{2}   & c_{6}  - c_{3}   \\
    \cdot & \cdot & c_{1} & c_{2} & c_{3} \\
    \cdot & \cdot & \cdot & c_{4} & c_{5} \\
    \cdot & \cdot & \cdot & \cdot & c_{6} \\
  \end{bmatrix}
$
}

\parbox{.06\linewidth}{\stepcounter{mac}(\themac)}%
\resizebox{.94\linewidth}{!}{%
$
  \begin{bmatrix}
     0 & -1 & -1 &  0 &  1 \\
     1 &  0 &  0 &  1 & -1 \\
     1 &  0 &  0 &  1 & -1 \\
     0 & -1 & -1 &  0 &  1 \\
    -1 &  1 &  1 & -1 &  0 \\
  \end{bmatrix}
\begin{bmatrix} v_{3}  - v_{2}   \\ v_{3}  - v_{1}   \\ v_{1} \\ v_{2} \\ v_{3} \end{bmatrix}
  \begin{bmatrix}
    c_{6}  -2  c_{5}   + c_{4}  &  - c_{3}   + c_{2}  + c_{6}  - c_{5}   & c_{3}  - c_{2}   & c_{5}  - c_{4}   & c_{6}  - c_{5}   \\
    \cdot & c_{6}  -2  c_{3}   + c_{1}  & c_{3}  - c_{1}   & c_{5}  - c_{2}   & c_{6}  - c_{3}   \\
    \cdot & \cdot & c_{1} & c_{2} & c_{3} \\
    \cdot & \cdot & \cdot & c_{4} & c_{5} \\
    \cdot & \cdot & \cdot & \cdot & c_{6} \\
  \end{bmatrix}
$
}

\parbox{.06\linewidth}{\stepcounter{mac}(\themac)}%
\resizebox{.94\linewidth}{!}{%
$
  \begin{bmatrix}
     0 & -1 & -1 &  1 &  1 \\
     1 &  0 & -1 &  0 &  0 \\
     1 &  1 &  0 & -1 & -1 \\
    -1 &  0 &  1 &  0 &  0 \\
    -1 &  0 &  1 &  0 &  0 \\
  \end{bmatrix}
\begin{bmatrix} v_{1} \\ v_{3}  + v_{2}  - v_{1}   \\ v_{1} \\ v_{2} \\ v_{3} \end{bmatrix}
  \begin{bmatrix}
    c_{1} & c_{3}  + c_{2}  - c_{1}   & c_{1} & c_{2} & c_{3} \\
    \cdot & c_{6}  +  2  c_{5}   -2  c_{3}   + c_{4}  -2  c_{2}   + c_{1}  & c_{3}  + c_{2}  - c_{1}   & c_{5}  + c_{4}  - c_{2}   & c_{6}  + c_{5}  - c_{3}   \\
    \cdot & \cdot & c_{1} & c_{2} & c_{3} \\
    \cdot & \cdot & \cdot & c_{4} & c_{5} \\
    \cdot & \cdot & \cdot & \cdot & c_{6} \\
  \end{bmatrix}
$
}

\bigskip 

\section{Periodicity shadows of size n=6}\label{sec:5} 

In the last section we present the longest list of all essential shadows of size $n=6$. Here we have $223$ essential 
shadows from the total number of $516$ (basic) tame periodicity shadows in $\bS(6)$. Note that there are 
$|\mathcal{S}(6)|=1290$ basic shades $6\times 6$. For the remaining shadows which are not essential, and the 
shades which are not shadows we refer to the complete lists in the Appendix. \medskip 

\setcounter{mac}{0}

\parbox{.07\linewidth}{\stepcounter{mac}(\themac)}%
\resizebox{.42\linewidth}{!}{%
$

$
}

\end{document}